\documentclass{gtmon_a}
\pdfoutput=1

\usepackage[all]{xy}

\proceedingstitle{The Zieschang Gedenkschrift}
\conferencestart{5 September 2007}
\conferenceend{8 September 2007}
\conferencename{Conference in honour of Heiner Zieschang}
\conferencelocation{Toulouse, France}

\editor{Michel Boileau}
\givenname{Michel}
\surname{Boileau}

\editor{Martin Scharlemann}
\givenname{Martin}
\surname{Scharlemann}

\editor{Richard Weidmann}
\givenname{Richard}
\surname{Weidmann}

\title{A classification of special $2$--fold coverings}

\author[A\,Bauval, D\,L\,Gon\c calves, C\,Hayat and M\,H\,P\,L\,Mello]{Anne
    Bauval\newline Daciberg L Gon\c calves\newline Claude
    Hayat\newline Maria Herm\'inia de Paula Leite  Mello}
\givenname{Anne}
\surname{Bauval}
\address{{\rm AB, CH:}\qua D\'epartement de Math\'ematiques\\
Laboratoire Emile Picard, UMR 5580\\\newline
Universit\'e Toulouse III\\
118 Route de Narbonne\\31400 Toulouse\\France\vspace{3pt}
\\\newline
{\rm DLG:}\qua Departamento de Matem\'atica - IME-USP\\\newline
Caixa Postal 66281\\   Ag\^encia Cidade de S\~ao Paulo\\
05311-970 - S\~ao Paulo - SP\\Brasil\vspace{3pt}\\\newline
{\rm MHPLM:}\qua Departamento de An\'alise Matem\'atica\\
Universidade Estadual do Rio de Janeiro\\\newline
Rio de Janeiro\\Brasil}
\email{bauval@picard.ups-tlse.fr}
\urladdr{}

\givenname{Daciberg L}
\surname{Gon\c calves}
\email{dlgoncal@ime.usp.br}
\urladdr{}

\givenname{Claude}
\surname{Hayat}
\email{hayat@picard.ups-tlse.fr}
\urladdr{}

\givenname{Maria Herm\'inia de Paula Leite}
\surname{Mello}
\email{mhplmello@ime.uerj.br}
\urladdr{}

\volumenumber{14}
\issuenumber{}
\publicationyear{2008}
\papernumber{02}
\startpage{27}
\endpage{47}

\doi{}
\MR{}
\Zbl{}

\arxivreference{} 


\received{7 June 2006}
\revised{10 Februray 2007} 
\accepted{16 January 2007} 
\published{29 April 2008}
\publishedonline{29 April 2008}
\proposed{}
\seconded{}
\corresponding{}
\version{}

\makeatletter
\def\cnewtheorem#1[#2]#3{\newtheorem{#1}{#3}
\expandafter\let\csname c@#1\endcsname\c@thm}

\makeatother  

\let\xysavmatrix\xymatrix
\def\xymatrix{\disablesubscriptcorrection\xysavmatrix}
\AtBeginDocument{\let\bar\wbar\let\tilde\wtilde\let\hat\what}
\AtBeginDocument{\renewcommand{\to}{\ensuremath{\longrightarrow}}}
\newcommand{\smrarrow}[1]{\smash{\stackrel{#1}{\smash{\longrightarrow}\vphantom{-}}}}

\numberwithin{equation}{section}
\newcommand{\intD}{\smash{\mskip4mu\mathring{\mskip-4mu\vrule width0pt height7pt depth0pt\smash{D}}}}
\AtBeginDocument{\renewcommand{\mod}{~\mathrm{mod}~}}

\newtheorem{thm}{Theorem}
\cnewtheorem{lem}[thm]{Lemma}
\cnewtheorem{prop}[thm]{Proposition}
\cnewtheorem{cor}[thm]{Corollary}

\theoremstyle{definition}
\cnewtheorem{defi}[thm]{Definition}
\cnewtheorem{nota}[thm]{Notation}
\cnewtheorem{rem}[thm]{Remark}

\newcommand{\Spin}{\rm Spin}

\begin{document}

\begin{asciiabstract} Starting with an $O(2)-principal fibration
over a closed oriented surface F_g, g>=1, a $2$--fold covering of the
total space is said to be \textit{special} when the monodromy sends
the fiber SO(2)\sim S^1 to the nontrivial element of Z_2. Adapting
D\,Jonhson's method [Spin structures and quadratic forms on surfaces,
J London Math Soc, 22 (1980) 365-373] we define an action of
Sp(Z_2,2g), the group of symplectic isomorphisms of (H_1(F_g;Z_2),.),
on the set of special $2$--fold coverings which has two orbits, one
with 2^{g-1}(2^g+1) elements and one with 2^{g-1}(2^g-1)
elements. These two orbits are obtained by considering Arf-invariants
and some congruence of the derived matrices coming from Fox Calculus.
Sp(Z_2,2g) is described as the union of conjugacy classes of two
subgroups, each of them fixing a special $2$--fold
covering. Generators of these two subgroups are made explicit.
\end{asciiabstract}

\begin{htmlabstract}
Starting with an SO(2)&ndash;principal fibration over a closed
oriented surface F<sub>g</sub>, g&ge; 1, a 2&ndash;fold covering of
the total space is said to be <em>special</em> when the monodromy
sends the fiber SO(2)&sim; S<sup>1</sup> to the nontrivial element of
<b>Z</b><sub>2</sub>. Adapting D&nbsp;Johnson's method [Spin
structures and quadratic forms on surfaces, J London Math Soc, 22
(1980) 365-373] we define an action of Sp(<b>Z</b><sub>2</sub>,2g),
the group of symplectic isomorphisms of
(H<sub>1</sub>(F<sub>g</sub>;<b>Z</b><sub>2</sub>),.), on the set of
special 2&ndash;fold coverings which has two orbits, one with
2<sup>g-1</sup>(2<sup>g</sup>+1) elements and one with
2<sup>g-1</sup>(2<sup>g</sup>-1) elements. These two orbits are
obtained by considering Arf-invariants and some congruence of the
derived matrices coming from Fox Calculus.
Sp(<b>Z</b><sub>2</sub>,2g) is described as the union of conjugacy
classes of two subgroups, each of them fixing a special 2&ndash;fold
covering. Generators of these two subgroups are made explicit.
\end{htmlabstract}

\begin{abstract} Starting with an $SO(2)$--principal fibration
over a closed oriented surface $F_g$, $g\geq 1$, a $2$--fold covering of
the total space is said to be \textit{special} when the monodromy sends the fiber
$SO(2)\sim S^1$  to the nontrivial element of $\mathbb{Z}_2$. Adapting D\,Johnson's
method \cite{Joh2}, we define an action  of $ Sp(\mathbb{Z}_2,2g)$, the group of
symplectic isomorphisms of $(H_1(F_g;\mathbb{Z}_2),.)$, on the set of special
$2$--fold coverings which has two orbits, one with $2^{g-1}(2^g+1)$
elements and one with $2^{g-1}(2^g-1)$ elements. These two orbits are obtained by
considering Arf-invariants and  some congruence of the derived matrices
coming from Fox Calculus.   $Sp(\mathbb{Z}_2,2g)$ is described as the union of
conjugacy classes of two subgroups, each of them  fixing a special
$2$--fold covering. Generators of these two subgroups are made explicit.
\end{abstract}

\begin{webabstract} Starting with an $SO(2)$--principal fibration
over a closed oriented surface $F_g$, $g\geq 1$, a $2$--fold covering
of the total space is said to be \textit{special} when the monodromy
sends the fiber $SO(2)\sim S^1$ to the nontrivial element of
$\mathbb{Z}_2$. Adapting D\,Johnson's method [Spin structures and
quadratic forms on surfaces, J London Math Soc, 22 (1980) 365-373] we
define an action of $ Sp(\mathbb{Z}_2,2g)$, the group of symplectic
isomorphisms of $(H_1(F_g;\mathbb{Z}_2),.)$, on the set of special
$2$--fold coverings which has two orbits, one with $2^{g-1}(2^g+1)$
elements and one with $2^{g-1}(2^g-1)$ elements. These two orbits are
obtained by considering Arf-invariants and some congruence of the
derived matrices coming from Fox Calculus.  $Sp(\mathbb{Z}_2,2g)$ is
described as the union of conjugacy classes of two subgroups, each of
them fixing a special $2$--fold covering. Generators of these two
subgroups are made explicit.
\end{webabstract}

\maketitle

\section{Introduction}\label{sec:intro}
We consider an  $SO(2)$--principal bundle over a closed oriented
surface $F_g$ of genus $g\geq 1$ as a $S^1$--principal bundle: $S^1\hookrightarrow
P\smrarrow{p}  F_g$.  A $2$--fold covering 
$\pi_{\varphi}\co E_\varphi \to P$ is said to be \textit{special} if  its
monodromy  $\varphi\co \pi_1 P\to \Z_2$ has the property that
$\varphi(u_0)=1$, where $1$ is   the nontrivial element of $\Z_2$, and  $u_0$ 
is the
image of the generator of $\pi_1 S^1$. The set ${\cal E}(q)=\{ \varphi\co
\pi_1 P\to \Z_2,\mid \varphi(u_0)=1\}$ is not empty if
and only if $q$, the Chern class   of the principal bundle, is even. This condition coincides with the vanishing of the second Stiefel--Whitney class of the $S^1$--principal bundle $S^1\hookrightarrow
P\smrarrow{p}  F_g$. In the sequel, it will be a running hypothesis that $q$ is even.
In \cite{DCHM} we obtained a presentation of $\pi_1E_\varphi$. For all $\varphi\in {\cal E}(q)$, these spaces $E_\varphi$ are isomorphic to the
total space of a $S^1$--bundle over $F_g$ classified by $\unfrac{q}{2}$.
The images $\pi_\varphi(\pi_1E_\varphi)$ are not conjugate subgroups of
$\pi_1 P$.  Nevertheless,  any $\varphi, \varphi'\in {\cal E}(q)$, $E_\varphi\to P$ and $E_{\varphi'}\to P$  are weakly equivalent in the sense that there
exists an automorphism $f$ of $\pi_1P$ such that $\varphi= \varphi'\circ
f$ (see \fullref{prop:we}).

The purpose of this work is to introduce on ${\cal E}(q)$ a supplementary structure obtained by an action of the symplectic group $Sp(H_1(F_g\Z_2),.)$.
The following theorem synthesizes the results obtained in \fullref{thm:to} and \fullref{thm:earf}.

\begin{thm}\label{thm:A} Let $\xi $ be a $S^1$--principal bundle over a closed surface 
$F_g$ of genus $g\geq 1$ with even Chern class $q$. Choosing a system of generators for $\pi_1F_g$ and $\pi_1P$ gives rise to a quadratic section  $s\co H_1(F_g;\Z_2)\to H_1(P;\Z_2)$ (see \fullref{prop:s}). The $s$--action of the symplectic group $Sp(H_1(F_g\Z_2),.)$ on the set ${\cal E}(q)$ of special $2$--fold coverings associated to the principal bundle $\xi$ produces two orbits, one with $2^{g-1}(2^g+1)$ elements and the other one with
$2^{g-1}(2^g-1)$ elements. The number of orbits and the number of elements in each orbit do not depend on the quadratic section~$s$.
\end{thm}

The quadratic section $s$ generalizes the work done by D\,Johnson \cite{Joh2} to any $S^1$--principal bundle over $F_g$ with even Chern class, when $\xi$ is associated to the tangent bundle of $F_g$, $g\geq 1$. Note that in this case the Chern class is always even.

One  motivation to study special $2$--fold coverings is that they can be
considered as $\Spin$--structures associated to an oriented $2$--vector
bundle over $F_g$ with even Chern class $q=2c$; see Milnor \cite{Mil} and the article by the last three authors \cite{DCM}. When this
oriented $2$--vector bundle is the tangent bundle and $F_g$ is orientable,
 Atiyah   \cite{A},  Birman and Craggs
\cite{BC} and Johnson \cite{Joh1,Joh3}  studied the
Torelli subgroup of the mapping class group of the surface $F_g$. In
these works, the  splitting of ${\cal E}(2c)$ into two classes is an important 
ingredient. Nevertheless the study of normal
fibrations defined by embeddings of a  surface in $\C^2$ (see Blanl{\oe}il and Saeki \cite{BS}) shows that
it is also worthwhile to start with any oriented  $S^1$--principal bundle
over  $F_g$.

The results of \fullref{thm:A} are obtained in two different ways. The quadratic section $s$ allows us to consider the set ${\cal E}(q)$ as the set of quadratic forms over $( H_1(F_g;\Z_2), .)$  where the symbol ``$.$'' is the intersection product. The associated Arf-invariant gives the counts of orbits and elements in each orbit. Considering the elements of ${\cal E}(q)$ as $2$--fold coverings leads us to use Crowell and Fox calculus and to define
congruence of the associated derived matrices (\fullref{defi:congr}). This congruence gives a classification of the $\Z_2[\Z_2]$--module structure of $H_1(E_\varphi,(E_\varphi)_0;\Z), \varphi\in {\cal E}(q)$ (\fullref{thm:fpsi}).

As shown by Atiyah \cite{A}, each symplectic automorphism fixes a quadratic form. In \fullref{cor:sp01}, $Sp(H_1(F_g\Z_2),.)$ is described as the union of conjugacy classes of two subgroups, each of them fixing a special $2$--fold covering. Generators of these two subgroups are made explicit in \fullref{thm:spo} and in \fullref{thm:sp1}.

\subsubsection*{Acknowledgements} This project started in 2002 and the last three authors were partially supported by the project CAPES-COFECUB number 364/01.

\section{First part}\label{sec:1P}

\subsection[Notation for the generators; introduction to special 2-fold coverings]{Notation for the generators; introduction to special $2$--fold
coverings}

  For $\pi_1(F_g,x)$ we take the usual presentation
$$\pi_1(F_g,x) = \Big\langle  x_1,\cdots,x_{2g} \Big|
              \prod_{j=1}^{g} [ x_{2j-1},x _{2j}]\Big\rangle .$$
In  $\pi_1(P,y)$ we choose elements
$\{ u_1,\cdots, u_{2g}\}$ such that
$p_\sharp( u_i)=x_i$. Let us fix ${\bf U}:=\{\{u_i\}_{1\leq i \leq
2g},u_0\}$ where $u_0$ is  a fixed  generator of the fiber of $p$.
The presentation of $\pi_1 (P,y)$ is:
$$\pi _1(P,y)=\Big\langle {\bf U}\Big| R_i=[u_i,u_0],1\leq i\leq 2g;
R_0=\prod_{\ell=1}^g[u_{2\ell -1},u_{2\ell }]u_0^{q}\Big\rangle.$$

\begin{defi}\label{defi:special} Let $u_0$ be the element of $\pi_1 P$
obtained from the fiber of $p$ and  consider the exact sequence
associated to a $2$--fold covering
$$1\to \pi_1E_{\varphi } \stackrel{\pi_\varphi}\longrightarrow \pi_1
P\stackrel{\varphi}\longrightarrow \Z_2\to 0.$$
When $\varphi(u_0)=1$, the nontrivial element of $\Z_2$, we will say
that the $2$--fold covering
$\pi_{\varphi}\co E_\varphi \to P$ is \textit{special}.
\end{defi}

There exists a one-to-one correspondence between the set of all special $2$--fold coverings
$\pi_{\varphi}\co E_\varphi \to P$   and 
the set ${\cal E}(q)=\{ \varphi\co \pi_1 P\to \Z_2 | \varphi(u_0)=1\}$, which 
corresponds bijectively to the set of
$\Spin$--structures associated to the oriented $2$--vector bundle over $F_g$
with  Chern class $q$. This set ${\cal E}(q)$ is not empty if and only if
$q$    is even (see the presentation of $\pi_1P$ given above). This
condition is valid throughout this work and  coincides with the vanishing of the second Stiefel--Whitney class of the $S^1$--principal bundle: $S^1\hookrightarrow
P\smrarrow{p}  F_g$. The set  ${\cal E}(q)$ has $2^{2g}$ elements.

One important property of these  special $2$--fold coverings  is that they
have isomorphic fundamental group \cite{DCHM}
$$\pi _1E_\varphi=\Big\langle y_1,\cdots,y_{2g},k\Big| [y_i,k],1\leq i\leq 2g;
\prod_{\ell=1}^g[y_{2\ell -1},y_{2\ell }]k^\frac{q}{2}\Big\rangle.$$
The injection  $\pi_\varphi\co \pi_1E_\varphi \to\pi_1 P$ is defined
by $\pi_\varphi(y_i)=u_i$,  if $\varphi(u_i)=0$, or
$\pi_\varphi(y_i)=u_iu_0^{-1}$, if $\varphi(u_i)=1$, and
$\pi_\varphi(k)=u_0^2$. There are $2^{2g}$ injections of this
type defining $2^{2g}$  images $\pi_\varphi(\pi_1E_\varphi)$ which are not
conjugate subgroups in $\pi_1P$. To be convinced of this fact, let us remark that
$\pi_\varphi(\pi_1E_\varphi)=\ker\varphi$, hence is a normal subgroup of $\pi_1P$, and 
 $\varphi=\varphi'$ if and only if 
$\ker\varphi=\ker\varphi'$.

\begin{rem}\label{rem:rem} Let us denote by $E^m$
the total space of a $S^1$--fibration over $F_g$ classified by the
integer $m$. Each $E^c$, with $c$ an odd integer, is the start of an
infinite graph with vertices  $E^{2^nc}$, and $2^{2g}$ arrows 
$E^{2^nc}\to E^{2^{n+1}c}$  the projections of nonisomorphic special
$2$--fold coverings.
\end{rem}

\begin{prop}\label{prop:1we} Two special $2$--fold coverings $E_\varphi \to
P$ and $E_{\varphi'} \to P$,  are always weakly
equivalent in the sense that there exists an automorphism $f$ of
$\pi_1P$ such that $\varphi= \varphi'\circ f$.
\end{prop}

Instead of proving this proposition, we will prove a stronger one,
\fullref{prop:we} in \fullref{arf}, where we impose 
$f$ to be a lift of an automorphism of $\pi_1F_g$. Let us recall some
facts about these different lifts.

\begin{lem}\label{lem:dehn}
{\rm (1)}\qua Let $\mathrm{Homeo}^+(F_g)$ be the group of homeomorphisms of $F_g$ preserving
the orientation. The projection
$$\mathrm{Homeo}^+(F_g)\to Sp(H_1(F_g;\Z_2),.)$$
is an epimorphism.

{\rm (2)}\qua An orientable homeomorphism of $F_g$ admits a lift as an orientable
fiber homeomorphism of $P$.
\end{lem}

\begin{proof}
\textbf{(1)}\qua The group of symplectic isomorphisms of  $(H_1(F_g;\Z_2),.)$ is generated by
the transvections, which are transformations  of the form
 $A(x)=x+(x.a)a$ for some vector $a$ \cite{omeara}.  These transvections
define Dehn twists, which are orientable homeomorphisms of the surface
$F_g$ \cite{ZVC}.

\textbf{(2)}\qua  Cutting the surface $F_g$ along a cut system produces a $4g$--polygon
$Y$. Let $D$ be a disk  in the interior of the polygon $Y$. The
restriction to $Y- \intD$ of the $S^1$--fibration $P$ is homeomorphic to $(Y-\intD)\times S^1$.
 To the  boundary of the hole, has to be
attached  a torus $D\times S^1$ after $q$ turns, where  $q$ is  the Chern 
class of $P$. Let $f$ be an orientable homeomorphism of $F_g$, we define $\hat
f$ to be $\smash{f|_{\smash{Y-\mathring{D}}}}\times \mathrm{id}$. The curve $f(\partial D)$ is a simple
closed curve. After turning $q$ times, the gluing of $\smash{{\hat f}}( (Y-\intD)\times
S^1)$ with $f(D)\times S^1$ is homeomorphic to $P$.
\end{proof}

The above results and considerations suggest that there exists an
action  of the group of symplectic isomorphism of $(H_1(F_g;\Z_2),.)$ on 
${\cal E}(q)$.

\subsection[Action of Sp(H1(Fg;Z2),.) on E(q)]{Action of $Sp(H_1(F_g;\Z_2),.)$ on ${\cal E}(q)$}

When $P$ is the $S^1$--principal bundle associated to the tangent bundle
of $F_g$,  Johnson defines an action of $Sp(H_1(F_g;\Z_2),.)$ which has two
orbits \cite{Joh2}. The definition of this action is given by means of a
choice of a section of the projection $H_1(P,\Z_2)\to
H_1(F_g,\Z_2)$. In  \cite{Joh2}, the section reflects the geometry of
the tangent bundle. We adapt this construction to make it  work for
any  oriented $S^1$--principal bundle over $F_g$.

\subsection[Johnson's lift of p-star: H1(P;Z2) to
H1(Fg;Z2)]{Johnson's lift of $p_\star\co  H_1(P;\Z_2)\to
H_1(F_g;\Z_2)$}\label{sec:s}

\begin{nota}\label{nota:hP}Let us denote by $h_M$ the composition $\pi_1M\rightarrow H_1(M;\Z)\rightarrow
H_1(M;\Z_2)$, where the first morphism is the
 Hurewicz epimorphism. An element $\varphi\in {\cal E}(q)$  determines a unique
$\tilde\varphi\co H_1(P;\Z_2)\to\Z_2$ such that
$\varphi=\tilde\varphi\circ h_P$ and $\varphi(u_0)=\tilde\varphi\circ
h_P(u_0)=1$, the nontrivial element of $\Z_2$. This allows us to
identify ${\cal E}(q)$ with $\{ \tilde\varphi\co H_1(P;\Z_2)\to  \Z_2\mid \tilde\varphi(u_0)=1\}$.
\end{nota}

The family $\boldsymbol{\sigma}=
\{ \sigma_i\}_{1\leq i\leq 2g}$, $\sigma_i:=h_{F_g}(x_i)$ where
$\{x_i\}$ are the fixed generators  of $\pi_1F_g$, is a symplectic basis
in $(H_1(F_g;\Z_2),.)$ where $.$ is the intersection product. The
family $\boldsymbol{\nu}:= \{\nu_i\}_{0\leq i\leq 2g}; \nu_i:=h_P(u_i)$ is a
basis of  $H_1(P;\Z_2)$.

\begin{prop}\label{prop:s} Choose a   family 
$\{s_i\}_{1\leq i \leq 2g}$  in  $\oplus_{0\leq i\leq 2g}\nu_i\Z_2=H_1(P;\Z_2)$ such that
$p_\star (s_i)=\sigma_i$ from the $2^{2g}$ possible choices. 
Then the following holds:

{\rm (1)}\qua $\{ \{s_i\}_{1\leq i \leq 2g},\nu_0\}$ is a  basis of
$H_1(P;\Z_2)$.

{\rm (2)}\qua For all $i$, $s_i=\nu_i+r_i\nu_0, r_i\in \Z_2$, so $2^{2g}$
possible choices  for $\{s_i, 1\leq i\leq 2g\}$.

{\rm (3)}\qua  There exists a
unique  map
 $$s\co \oplus_{1\leq i\leq 2g} \sigma_i\Z_2=H_1(F_g;\Z_2)\to\oplus_{0\leq i\leq 2g}\nu_i\Z_2=H_1(P;\Z_2),$$
 defined by $s(\sigma_i)=s_i, 1\leq i\leq 2g$ such that
for all $a,b\in H_1(F_g;\Z_2)$
 \begin{eqnarray}\label{eq:*}
 s(a+b)=s(a)+s(b)+(a.b)\nu_0.
\end{eqnarray}
\end{prop}

\begin{nota}\label{nota:qs} The map $s$ obtained in \fullref{prop:s} will be called a \textit{quadratic section}.
\end{nota}
 
\begin{proof}[Proof of \fullref{prop:s}]
\textbf{(1)}\qua If $\Sigma \alpha _is(\sigma_i)+\gamma \nu_0=0$, then
$p_\star(\Sigma \alpha _is(\sigma_i)+\gamma \nu_0)= \Sigma \alpha
_i\sigma_i=0$; so
$\alpha _i=0$ and $\gamma =0$.
This implies that $\{ \{s(\sigma_i)\}_{1\leq i \leq 2g},\nu_0\}$ is a 
basis of $H_1(P;\Z_2)$.

\textbf{(2)}\qua This is true because $\ker p_\star=\langle\nu_0\rangle$.

\textbf{(3)}\qua Take $a=\Sigma a_i\sigma_i,a_i\in \Z_2$. Because of condition
\eqref{eq:*}, we must 
 define $s(a)$ by:
\begin{eqnarray*}
s(a)=\Sigma a_is(\sigma_i)+ (\Sigma a_{2i-1}a_{2i})\nu_0.
\end{eqnarray*}
 Now, if the map $s$ is defined by this equation, then
$$\begin{array}{ll}
s(\Sigma a_i\sigma_i+\Sigma b_i\sigma_i) &= s(\Sigma (a_i+
b_i)\sigma_i) \\
                                         &=\Sigma (a_i+ b_i)s(\sigma_i)+
[(\Sigma (a_{2i-1}+ b_{2i-1})(a_{2i}+ b_{2i})]\nu_0 \\
                                         &= s(\Sigma a_i\sigma_i)+s(\Sigma
b_i\sigma_i)+[\Sigma (a_{2i}b_{2i-1}+ a_{2i-1}b_{2i})]\nu_0 \\
                                         & = s(\Sigma a_i\sigma_i)+s(\Sigma
b_i\sigma_i)+ (\Sigma a_i\sigma_i).(\Sigma b_i\sigma_i)\nu_0,
\end{array}
$$
 because the coefficients are in $\Z_2$.
\end{proof}

\begin{rem}\label{rem:ns}
In the case where the fiber bundle $P$ is the $S^1$--principal bundle
associated to the tangent bundle of the surface $F_g$, the geometry
imposes the choice of $r_i=1,1\leq i \leq 2g$  \cite{Joh2}.
Hence, if  necessary, it is possible to normalize the choice of the
maps $s$ imposing this condition on the family $r_i$ as in Arf \cite{A}.

Let U be a set of generators of $H_1(P;Z_2)$. We also impose that for each 
chosen generator of $\pi_1F_g$ there will be one element in $U$ which is a 
lift of it. Hence
 two such  systems of generators ${\bf U}$ and ${\bf U'}$ of  $\pi_1P$ are
related by $u'_i=u_0^{-\alpha_i}u_i$ (or equivalently
$u_i=u_0^{\alpha_i}u'_i$), $\alpha_i\in \{0,1\}$.  An element 
$\varphi\in {\cal E}(q)$ is then changed into
$\varphi'(u_i)=\varphi(u_i)+\alpha_i$ and $\varphi'(u_0)=\varphi(u_0)$.
Note that such a change of generators is equivalent to a change of the quadratic section $s$.
\end{rem}

\begin{defi}\label{defi:tildef} Let  $A$ be the symplectic
matrix $(a_{ij})_{i,j\leq 2g}$ written in the basis $\boldsymbol{\sigma}$, of
a symplectic
isomorphism $f\co H_1(F_g;\Z_2)\to H_1(F_g;\Z_2)$. We define
 $$ f_s\co \oplus_{0\leq i\leq 2g}\nu_i\Z_2=H_1(P;\Z_2)\to \oplus_{0\leq i\leq 2g}\nu_i\Z_2=H_1(P;\Z_2)$$
 by linearity
from
$$ f_s(s(\sigma_i)):=s(f(\sigma_i)),  f_s(\nu_0):=\nu_0.$$
 The matrix of $ f_s$ in the basis $\boldsymbol{\nu}$ is
$\big(\begin{smallmatrix}
                                  A&0\\
                                W&1
                                 \end{smallmatrix}\big)$
where $W$ is a line with $2g$ terms $w_j=\Sigma a_{ij}r_i+S_j+r_j$, $
S_j=\Sigma a_{2i,j}a_{2i-1,j}$,  $r_j\nu_0=s(\sigma_j)+\nu_j. $
\end{defi}

Notice that $f_s\circ s=s\circ f$. We have $(f_1f_2)_s=(f_1)_s(f_2)_s$
and $(\mathrm{id}_{Sp(\Z _2,2g)})_s=\mathrm{id}_{ Sl(\Z _2,2g+1)}$.
This proves the following proposition, where $Sp(\Z _2,2g)$  denotes the
group of the symplectic $2g\times 2g$ matrices with
coefficients in $\Z _2$:
\eject 
\begin{prop}\label{prop:J}The injective map
$$
\begin{array}{rcl}
J\co Sp(\Z _2,2g)&\to& Sl(\Z _2,2g+1)\\
A&\mapsto &\tilde A=\left( \begin{array}{cc}
                                  A&0\\
                                 W&1
                                 \end{array}
                                  \right)
\end{array}
$$
 with
$A=(a_{ij})_{i,j\leq 2g}$ and $W=(w_1\cdots w_{2g})$ where $w_j=\Sigma
a_{ij}r_i+S_j+r_j$ and $S_j=\Sigma a_{2i,j}a_{2i-1,j}$
is a monomorphism.
\end{prop}
 
\subsection[s-Relation between two special 2-fold coverings]{$s$--Relation between two special $2$--fold coverings}\label{sec:sr}

\begin{defi}\label{defi:rsp}
Two special $2$--fold coverings $\varphi$ and $\varphi'$ are $s$--related if
there exists a symplectic isomorphism
$f\co H_1(F_g;\Z_2)\to H_1(F_g;\Z_2)$ such that
$$\tilde\varphi=\tilde\varphi'\circ  f_s,$$
with $ f_s$ given in \fullref{defi:tildef}, or equivalently: $\tilde\varphi\circ s=\tilde\varphi'\circ s\circ f$.
\end{defi}

\begin{prop}\label{prop:relation}
Two special $2$--fold coverings $\varphi$ and $\varphi'$ are $s$--related if
and only if
  \begin{equation}\label{eqn:star}
\forall j , \varphi(u_j)=\Sigma_{i=1,\ldots, 2g}a_{ij}\varphi' (u_i)+ w_j,
\end{equation}
where  $A=(a_{ij})$ is  the matrix of a symplectic isomorphism in the basis $\sigma$, and $ w_j=\Sigma
a_{ij}r_i+\Sigma a_{2i,j}a_{2i-1,j}+r_j$, $r_j$ determined by the choice
of $s$.\qed
\end{prop}

Let us recall or introduce some terminology needed for \fullref{thm:to}.
\begin{itemize}
\item For any automorphism $K$ of $H_1(F_g;\Z_2)$, a \textit{lift} of $K$ is an automorphism $k$ of $\pi_1P$ such that $h_P\circ k=K\circ h_P$, where $h_P$ is defined in \fullref{nota:hP}.
\item  Two special $2$--fold coverings $\pi_\varphi\co E_\varphi\to P$ and $\pi_{\varphi'}\co E_{\varphi'}\to P$ are \textit{weakly equivalent} by
$k\in \mathrm{Aut}(\pi_1(P))$ if and only if
$\varphi=\varphi'\circ k$.
\end{itemize}

\begin{thm}\label{thm:to}
{\rm (1)}\qua For any symplectic automorphism $f$ of $H_1(F_g;\Z_2)$, there exists a lift $f_\sharp\co \pi_1P\to \pi_1P$ of
$f_s$ (\fullref{defi:tildef}).

{\rm (2)}\qua For any such $f$ and $f_\sharp$, $\varphi,\varphi'\in {\cal E}(q)$ are $s$--related
by $f$ if and only if
$\pi_\varphi,\pi_{\varphi'}$ are weakly equivalent by $f_\sharp$.
\end{thm}

\begin{proof}
\textbf{(1)}\qua Using geometric arguments we proved in \fullref{lem:dehn} 
that
there exists an automorphism $g$ of $\pi_1P$ such that $g(u_0)=u_0$ and
$f\circ h_{F_g}\circ p_\sharp=h_{F_g}\circ p_\sharp\circ g$.  
This implies that the morphism $f_s\circ h_P-h_P\circ g\co\pi_1(P)\to
H_1(P;\Z_2)$
takes its values in the subgroup $\ker(p_*)=\Z_2\nu_0$, ie, it is of the
form
$x\mapsto\bar\rho(x)\nu_0$ for some morphism $\bar\rho
\co\pi_1(P)\to\Z_2$.
If we are able to construct a lift $\rho \co\pi_1(P)\to\Z$ of $\bar\rho$,
then we just have to define
$f_\sharp\co\pi_1(P)\to\pi_1(P)$ by $f_\sharp(x)=g(xu_0^{\rho(x)})$ to get an
automorphism $f_\sharp$ of $\pi_1(P)$
satisfying $h_P\circ f_\sharp=f_s\circ h_P$.
In order to construct such a lift $\rho$, notice that $\bar\rho$
factorizes through $H_1(P;\Z)$:
denoting by $hz\co\pi_1(P)\to H_1(P;\Z)$ the Hurewicz morphism, we have
$\bar\rho=\bar r\circ hz$ for some morphism $\bar
r\co H_1(P;\Z)\to\Z_2$,
for which we want a lift $r\co H_1(P;\Z)\to\Z$.
There are many such $r$'s, since $\bar r(hz(u_0))=\bar\rho(u_0)=0$
and
$H_1(P;\Z)/\langle hz(u_0)\rangle\simeq H_1(F_g;\Z)$ is a free $\Z$--module.

\textbf{(2)}\qua  Recall that $h_P$ is an epimorphism, hence we have the following equivalences:
$$\tilde\varphi'\circ
f_s=\tilde\varphi\Leftrightarrow\tilde\varphi'\circ f_s\circ
h_P=\tilde\varphi\circ h_P\Leftrightarrow\tilde\varphi'\circ
h_P\circ f_{\sharp}=\tilde\varphi\circ h_P\Leftrightarrow\varphi'\circ
f_{\sharp}=\varphi.\proved$$
\end{proof}

\subsection{Arf type invariant}\label{arf}

The purpose  of this section is to prove that there are two orbits
under the $s$--action, one with $2^{g-1}(2^g+1)$ elements and one with
$2^{g-1}(2^g-1)$ elements. We use the quadratic section $s$ defined and fixed  in
the above subsection  to associate bijectively a special $2$--fold covering
$\varphi$ and a quadratic form $\omega_\varphi=\tilde\varphi\circ s$.

Let $\{\sigma_1,\sigma_2,\cdots,\sigma_{2g-1},\sigma_{2g}\}$ be a
symplectic basis
of $(\Z_2^{2g},.)$. This means that
$\sigma_{2i-1}\,{.}\,\sigma_{2i}=\sigma_{2i}\,{.}\,\sigma_{2i-1}=1, 1\leq i\leq g$,
and all the
others $\sigma_i.\sigma_j=0$. The Arf-invariant of a quadratic form
$\omega\co
 (\Z_2^{2g},.)\to
\Z_2$  is defined by
$$\alpha(\omega)=\Sigma\omega(\sigma_{2j-1})\omega(\sigma_{ 2j}).$$

\begin{thm}\label{thm:earf}  Two special $2$--fold coverings $\varphi$ and
$\varphi' $ are $s$--related  if and only if the Arf-invariants of
$\omega_\varphi$ and $\omega_{\varphi'}$  are equal \cite{Arf},
explicitly:
$$\Sigma{\tilde\varphi}(s(\sigma_{2j-1})){\tilde\varphi}(s(\sigma_{2j}))
=\Sigma{\tilde\varphi'}(s(\sigma_{2j-1})){\tilde\varphi'}(s(\sigma_{2j})).$$
\end{thm}

\begin{proof}
\fullref{prop:s} proved that $s(a+b)=s(a)+s(b)+(a.b)\nu_0$, so
$\varphi$ determines a quadratic form
$$
\begin{array}{rcl}
\omega_\varphi\co \Z _2^{2g}&\to& \Z _2\\
a&\mapsto & \omega_\varphi(a)={\tilde\varphi}(s(a)).
\end{array}
$$
A quadratic form $\omega$ determines $\varphi\in {\cal E}(q)$ by
$\tilde\varphi(s(\sigma_i))=\omega(\sigma_i)$ and
$\tilde\varphi(\nu_0)=\nu_0$. Two special $2$--fold coverings $\varphi$ and
$\varphi' $ are $s$--related (\fullref{defi:rsp}) if and only if
there exists  a symplectic map
$f\co H_1(F_g,\Z_2)\to H_1(F_g,\Z_2)$ such that 
$\omega_\varphi=\omega_{\varphi'}\circ f$, which is equivalent to the
equality of the Arf-invariants of $\omega_\varphi$ and
$\omega_{\varphi'}$ \cite{Arf}. We give below a short proof of this
classical property.
\end{proof}

 \begin{prop}\label{prop:aa}
There exists a symplectic map $f\co (\Z_2^{2g},.)\to(\Z_2^{2g},.) $
such that
$\omega=\omega'\circ f$  if and only if $
\alpha (\omega)=\alpha (\omega')$. We will denote this by
$\omega\sim\omega'$.
\end{prop}

\begin{proof}
Let $\omega ,\omega '\co(\Z_2^{2g},.)\to\Z_2$ be any two quadratic forms.
Their
difference is a linear form
$$\omega '(x)-\omega (x)=V.x.$$
By an elementary computation we have:
$$\alpha (\omega')-\alpha(\omega)=\omega(V).$$
For any vector $Y$, let us denote by $T_Y$ the symplectic transvection defined by
$T_Y(x)=x+(Y.x)Y$. We then obtain
$\omega(T_Y(x))=\omega (x)+\omega((Y.x)Y)+Y.x$, hence
$$(\omega\circ T_Y)(x)-\omega(x)=(1+\omega(Y))Y.x.$$
Using these two equations we deduce:
\begin{itemize}
\item $\alpha(\omega')=\alpha(\omega)\Rightarrow\omega(V)=0\Rightarrow
\omega \circ T_V-\omega =V.\ =\omega '-\omega \Rightarrow\omega \circ
T_V=\omega'\Rightarrow \omega'\sim\omega$.
\item Conversely, $\omega '=\omega \circ T_Y\Rightarrow V=(1+\omega
(Y))Y\Rightarrow\alpha (\omega ')-\alpha (\omega )=(1+\omega (Y))\omega
(Y)=0$.  Hence (since transvections generate the group of symplectic
isomorphisms \cite{omeara}) $\omega '\sim\omega \Rightarrow \alpha
(\omega ')=\alpha
(\omega )$.\proved
\end{itemize}
 \end{proof}

The following proposition will prove a stronger property than weak
equivalence  for any pair of special $2$--fold coverings:

\begin{prop}\label{prop:we} Given two special $2$--fold coverings 
$E_\varphi\to P, E_{\varphi'}\to P$,  it is possible to choose a quadratic section
$s(\varphi,\varphi')$ such that these $2$--fold coverings are 
$s(\varphi,\varphi')$--related (\fullref{defi:rsp}).
\end{prop}
\begin{proof}
First, it is possible to choose   a quadratic section
 $s=s(\varphi,\varphi')$ such that
$\alpha (\tilde{\varphi}\circ s)=0=\alpha (\tilde{\varphi}'\circ s)$.
In fact, because
$\alpha (\tilde{\varphi}\circ
s)=\smash{\sum_{i=1}^g}(\tilde{\varphi}(\nu_{2i-1})+r_{2i-1})(\tilde{\varphi}(\nu_{2i})
+r_{2i})$
(the same for $\varphi '$), it is enough to choose for example
$r_i=\tilde{\varphi}(\nu_i)$ for $i$ odd and
$r_i=\tilde{\varphi}'(\nu_i)$ for $i$ even. By \fullref{prop:aa}  or \cite{Arf} there
exists
$f\in Sp(H_1(F_g;\Z_2),.)$ such that
$\tilde\varphi \circ s=\tilde\varphi'\circ  s\circ f$. 
\end{proof}

\subsection[Subgroups Sp-omega (Z2,2g) of the symplectic automorphisms which fix a quadratic form omega]{Subgroups $Sp_\omega (\Z_2,2g)$ of the symplectic
automorphisms which fix a quadratic form $\omega$}\label{sec:spomega}

\subsubsection{Generators}

As shown by Atiyah \cite{A}, each symplectic
automorphism fixes a quadratic form  $\omega $.
Let us study the subgroup $Sp_\omega (\Z_2,2g)$ of symplectic
automorphisms which fix $\omega $ (this $\omega $ may be of the form
$\omega _\varphi :=\tilde\varphi \circ s$).

It suffices to study the two subgroups $Sp_i$ ($i=0$ or $1$)
corresponding to $\omega _i$, with $\omega _0(x):=\sum x_{2k-1}x_{2k}$
and
$\omega _1(x):=\omega _0(x)+x_1+x_2$. Then, if $\alpha (\omega
)=i$, $Sp_\omega $ is a conjugate of $Sp_i$ (by any $f\in Sp(\Z_2,2g)$
such that
$\omega=\omega _i\circ f$).

\begin{lem}\label{lem:15} The actions of $Sp_0$ on $H_0:=\{x\neq
0,\omega_0(x)=0\}$ and on
$H_1:=\{x,\omega _0(x)=1\}$ are transitive.
\end{lem}

\begin{proof}  We assume that $g>1$ ($g=1$ is obvious). Note that $Sp_0$
contains all symplectic permutations, and all transvections $T_u$ such
that
$\omega _0(u)=1$.

If $x\in H_0$, since $x\neq 0$, up to
some symplectic permutation, we may assume that $x.e_1=1$. Let $u:=x+e_1$. Then $\omega
_0(u)=1$ and $T_u(e_1)=x$.

 If $x\in H_1$, we have:

Case 1:\qua
If $x.(e_{2k-1}+e_{2k})=1$ for some $k$, up to
some symplectic permutation, we may assume that $k=1$. Let $u:=x+e_1+e_2$.
Then $\omega _0(u)=1$ and $T_u(e_1+e_2)=x$.

Case 2:\qua If $x.(e_{2k-1}+e_{2k})=0$ for all $k$'s. Since $x\neq 0$,
up to some symplectic permutation, we may assume that $x.e_1=1$. Let
$u':=e_1+e_3+e_4$ (hence $\omega _0(u')=1$) and $x':=T_{u'}(x)=x+u'$.
Then
$x'.(e_1+e_2)=u'.(e_1+e_2)=1$ hence we are led to the first case.
\end{proof}

\begin{thm}\label{thm:spo} Any element of $Sp_0$ is a product of:
\begin{itemize}
\item[\rm(1)] symplectic permutations,
\item[\rm(2)] (if $g\geq 2$) the matrix $B_1:=\left( \begin{array}{cc} A_1&0\\
0&I_{2g-4}\end{array}\right)$
with $A_1=\left( \begin{array}{cccc}1&0&0&0\\ 0&1&0&1\\ 1&0&1&0\cr
0&0&0&1\end{array}\right)$.
\end{itemize}
\end{thm}

\begin{proof} Take $g>1$ ($g=1$ is obvious)  and assume the property
true
for $g-1$. Call ``type R'' all matrices of the form $\big(
\begin{smallmatrix}I_2&0\\
0&A\end{smallmatrix}\big)$ (which, by induction hypothesis, are products of
these
generators).
Let $\gamma \in Sp_0$ and $V:=\mathrm{Vect}(e_1,e_2)$.

Case 0:\qua $\gamma (V)=V$.\qua Then $\gamma $ fixes or exchanges $e_1$ and
$e_2$; hence (up to some product by a symplectic transposition) $\gamma $
is of type $R$.

Case 1:\qua $\gamma (V)\neq V$ but there exists a (nonzero) $x\in
V$ such that $\gamma (x)\in V$.
\begin{enumerate}
\item[1.1:] $x=e_1$ or $e_2$.\qua Up to symplectic
permutation(s), $\gamma (e_2)=e_2$. Then $\gamma(e_1)=y+z$ with $y\in
V,z\in V^\perp$, $z\neq 0$, $y.e_2=1$ (hence $y=e_1$ or $e_1+e_2$), and
$\omega _0(z)=\omega _0(y)$.
\begin{itemize}
\item[1.1.1:] $y=e_1$, $\omega_0(z)=0$.\qua Hence by the lemma we may
assume $z=e_3$ (up to some product by a type R matrix). In this case,
$B_1^{-1}\gamma $ is of type R.
\item[1.1.2:] $y=e_1+e_2$, $\omega
_0(z)=1$.\qua Hence (by the lemma again) we may assume $z=e_3+e_4$. In this
case,
$B_1^{-1}\gamma$ fixes $e_2$ and sends $e_1$ to $e_1+e_4$, hence it falls into
the subcase 1.1.1.
\end{itemize}
\item[1.2:] $x=e_1+e_2$.\qua For $i=1,2$, $\gamma (e_i)=y_i+z$
with
$y_i\in V,z\in V^\perp, z\neq 0, y_1+y_2=e_1+e_2, y_1.y_2=1$. Hence
(up to symplectic transposition) $y_i=e_i$, so that $\omega _0(z)=0$,
hence,
by the lemma again, we may assume that  $z=e_3$. In that case, $B_1^{-1}\gamma$ fixes
$e_1$;  hence it belongs to the subcase 1.1 (or to case 0).
\end{enumerate}

Case 2:\qua $\gamma (x)\in
V^\perp$ for some nonzero $x\in V$.\qua By the lemma we may assume (up to
some
product by a type R matrix) that $\gamma (x)=e_3$ or $\gamma
(x)=e_3+e_4$
(depending whether $\omega _0(x)$ equals $0$ or $1$). By symplectic
permutation the situation is  reduced to case 0 or 1.

Case 3:\qua None of the three nonzero elements of $V$ is sent by $\gamma
$ to $V\cup V^\perp$.\qua Let $\gamma(e_1)=y+z,\gamma (e_2)=y'+t$ with
$y,y'\in V,z,t\in V^\perp$. Then $y,y'$ are nonzero and distinct, hence
at least one of them equals some $e_i$ (with $i=1$ or $2$). We may assume that
$\gamma (e_1)=e_1+z$, hence $\omega _0(z)=0$. Since $z\neq 0$, we may
assume $z=e_3$.
Then $B_1^{-1}\gamma$ fixes $e_1$, hence it belongs to case 0 or 1.
\end{proof}

\begin{rem} A classical set of generators for the whole group $Sp(\Z_2,2g)$
consists of these generators of the subgroup $Sp_0$, and the matrix
$B_0$
corresponding to $A_0:=\big( \begin{smallmatrix}1&1\\
0&1\end{smallmatrix}\big)$ (cf O'Meara \cite{omeara}).
\end{rem}

\begin{thm}\label{thm:sp1} Any element of $Sp_1$ is a product of:
\begin{itemize}
\item[\rm(1)] elements of the subgroup $Sp(\Z_2,2)\times Sp_0(\Z_2,2g-2)$,
\item[\rm(2)] (if $g\geq 2$) the matrix $B_2:=\left( \begin{array}{cc}A_2&0\\
0&I_{2g-4}\end{array}\right)$
with $A_2=\left( \begin{array}{cccc}1&0&0&0\\ 0&1&0&1\\
1&0&1&1\\0&0&0&1\end{array}\right)$.
\end{itemize}
\end{thm}

\begin{proof} If $g=1$, $Sp_1=Sp(\Z_2,2)$.

Let $\gamma \in Sp_1$ and $V=\mathrm{Vect}(e_1,e_2)$.

Case 0:\qua $\gamma (V)=V$.\qua Then
$\gamma \in Sp(\Z_2,2)\times Sp_0$.

Case 1:\qua $\gamma (V)\neq V$ but there exists a (nonzero) $x\in
V$ such that $\gamma (x)\in V$.\qua Assume (up to products by elements of
$Sp(\Z_2,2)=GL(2,\Z_2)$) $\gamma (e_2)=e_2$ and $\gamma (e_1)=e_1+z$,
with
$z\in V^\perp$, nonzero, and such that $\omega _0(z)=0$. Assume moreover
(up to some product by an element of $Sp_0$, by \fullref{lem:15})
$z=e_3$. Then $B_2^{-1}\gamma $ belongs to the subgroup
$Sp_0(\Z_2,2g-2)$.

Case 2:\qua For some $x\in V$, $\gamma (x)\notin V\cup V^\perp$.\qua 
Using
the same arguments as above, we may assume that $\gamma (e_1)=e_1+e_3=B_2(e_1)$, hence
$B_2^{-1}\gamma$ satisfies the condition in case  0 or 1.

Case 3:\qua For all $x\in V,\gamma (x)\in V^\perp$.\qua We may assume that $\gamma
(e_1)=e_3+e_4=B_2(e_2+e_4)$, hence $B_2^{-1}\gamma$ satisfies the condition in 
case  2.
\end{proof}

For each $\omega$ such that $\alpha(\omega)=\alpha(\omega_i), i=0,1$, let us 
choose the transvection $T_{Y_\omega}$ where $Y_\omega$ is the vector
such that for all $x$, $\omega(x)-\omega_i(x)=Y_\omega.x$. Recall that we have shown in the proof
of \fullref{prop:aa} that
$\alpha(\omega)-\alpha(\omega_i)=\omega_i(Y_\omega)$. Now let us define the two
subsets  $\alpha_i:=\{Y\mid\omega_i(Y)=0\}. $ The
family  $\alpha_0$ has $2^{g-1}(2^g+1)$ elements and  $\alpha_1$ has
$2^{g-1}(2^g-1)$ elements. We get the corollary:

\begin{cor}\label{cor:sp01}
$$Sp(\Z_2,2g)= \bigcup_{Y\in \alpha_0}[T_Y^{-1}Sp_0
T_Y]\cup\bigcup_{Y\in \alpha_1}[T_Y^{-1} Sp_1T_Y]. $$
\end{cor}

The generators of $Sp_0$ and $Sp_1$ (Theorems \ref{thm:spo},
\ref{thm:sp1}) admit lifts, described for example in Zieschang, Vogt and Coldewey \cite{ZVC}, as
homeomorphisms of the surface $F_g$. When a quadratic section $s$ is chosen, we may
view these homeomorphisms as homeomorphisms fixing a
$\Spin$--structure associated to an oriented $2$--vector bundle over $F_g$
with Chern class equal to $q$.

\begin{cor}\label{cor:sp013} {\rm (1)}\qua Under the action defined in 
\fullref{defi:rsp} the set ${\cal E}(q)$ of special $2$--fold
coverings is divided into two orbits: ${\cal E}(q)^0$ with
$2^{g-1}(2^g+1)$ elements and ${\cal E}(q)^1$  with $2^{g-1}(2^g-1)$.

{\rm (2)}\qua The stabilizers of an element of ${\cal E}(q)^i$ is a conjugate of
$Sp_i, i=0,1$.
\end{cor}

\begin{rem} Let us emphasize  that after a change of the generators of
$\pi_1P$, which are lifts of the fixed generators of $\pi_1F_g$, or
after a change in the choice of the quadratic section $s$ (see 
\fullref{prop:s} (2)), only the number of orbits of ${\cal E}(q)$ and the
number of elements in each orbit do not change.
\end{rem}

\section{Second part}

\subsection{Derived matrix}

In this section we apply to the special $2$--fold coverings the classical
tools of Fox derivatives. We will give a description of the $\Z [\Z
_2]$--module
structure of $H_1(E_{\varphi },(E_{\varphi })_0;\Z)$, using the
Reidemeister method as referred in \cite[Chapter 9]{BZ} (also \cite{CF}),  where $(E_{\varphi })_0$ is the
fiber with two
 elements above the base point of $P$.
The exact sequence of the pair $(E_{\varphi }, (E_{\varphi })_0$ is:
$$0\to H_1(E_{\varphi };\Z)\to H_1(E_{\varphi},(E_{\varphi })_0;\Z) \to
\Z[\Z _2] \to \Z
\to 0.$$  A notion of congruence is defined on the matrices. It leads to
the same relation between the data as the necessary relations  to be
$s$--related \eqref{eqn:star} or Arf related \fullref{thm:earf}. The last
step is to add a $\ast$--product on $H_1(E_{\varphi},(E_{\varphi
})_0;\Z_2)$ and to find a relation between   the $\Z [\Z
_2]$--module
structures of $H_1(E_{\varphi},(E_{\varphi
})_0;\Z_2)$ and $H_1(E_{\varphi'},(E_{\varphi' })_0;\Z_2)$ when
$\varphi$ and $\varphi'$ are $s$--related.

\subsubsection{Summary of Crowell and Fox calculus}\label{CF}

Let $\varphi\in {\cal E}(q)$. In the exact sequence of the homotopy groups of
the special  $2$--fold covering
$\pi\co E_{\varphi} \longrightarrow P$:
$$0\to \pi _1 (E_{\varphi },x)\stackrel{\pi_{\sharp} }\longrightarrow
\pi _1(P,y) \stackrel{\varphi }\longrightarrow\Z _2\to 0$$
 the group $\Z _2$ is the multiplicative group of deck transformation of the
covering. Writing the ring
 $\Z [\Z _2]= \unfrac{\Z [t]}{(1-t^2)}$, the homomorphism
 $\varphi \co \pi _1P \to \Z_2 $ extends to a group ring morphism
 $\Z [ \pi _1P]\to \Z[\Z_2]$, also denoted by  $\varphi $. This morphism
 verifies in particular $\varphi (u_0)=t,\varphi(1)=1$, ($u_0$ coming
from the fiber $S^1$) and $\varphi (0)= 0$.

\subsubsection{Explicit computations of the derived matrix}\label{sec:qo}

Recall that we make  choices such that the presentation of  $\pi
_1(P,y)$ is:
 $$\pi _1(P,y)=\Big\langle{\bf  U}\Big| R_i=[u_i,u_0],1\leq i\leq 2g;
R_0=\prod^g_1[u_{2\ell -1},u_{2\ell }]u_0^{2c}\Big\rangle.$$
 Let $M_2$ be the free $\Z [\Z _2]$--module generated by the set ${\bf
R}:=\{R_i,0\leq i\leq 2g\}$ and $M_1$ the free $\Z [\Z _2]$--module
generated by  ${\bf U}=\{u_i,0\leq i\leq 2g\}$. Using the Fox derivation 
$\partial R_i/ \partial u_j$, 
a $\Z [\Z _2]$--morphism $d_\varphi\co (M_2,{\bf R})\to (M_1,{\bf U})$
is defined  by
$d_\varphi R_i=\sum_j m_{ji}u_j$, where
$m_{ji}=\varphi q(\partial R_i/\partial u_j)$ and $q$ is the ring
morphism obtained from the group projection from the free group generated by
the set ${\bf U }$ to $\pi _1(P,y)$. So there is an exact sequence
of $\Z [\Z _2]$--modules:
$$( M_2, {\bf R})\stackrel{d_\varphi }\longrightarrow (M_1,{\bf U})\to
(M_1/ \mathrm{Im}\ d_\varphi, {\bf \bar{ U}})\to 0,$$
 where ${\bf\bar{ U}}:=\{\bar{u_i}\}_{0\leq i\leq 2g}$,
$\bar{u_i}$ class of
$u_i$ modulo $ \mathrm{Im} d_\varphi$.

The structure of  $\Z[\Z _2]$--module of $M_1/ \mathrm{Im}\  d_\varphi$ is
denoted  by $H_\varphi$.

 Let $u$ be an element of $\pi_1(P,y)$ and select a loop $\alpha \in u$.
By the path-lifting property of covering spaces, there exists a unique
path $\alpha'\co I\to E$ such that the projection of $\alpha'$ is
$\alpha$ and $\alpha'(0)=y$. Its relative homology class is denoted by
$\tilde u$.
From \cite[Chapter 9]{BZ}, \cite{CF}, we know that there exists a $\Z[\Z
_2]$--isomorphism
\begin{equation*}
H_\varphi\to H_1(E_{\varphi},(E_{\varphi})_0;\Z),
{\bar u_i}\mapsto {\tilde u_i}.
\end{equation*}
Up to this isomorphism, we have to study the $\Z[\Z _2]$--module
$H_\varphi$.

We introduce the notation $n=(n_1,\cdots ,n_{2g})$ where
$n_i=0$ if $\varphi (u_i)=1$ and $n_i=-1$ if $\varphi (u_i)=t$. For
convenience, we also denote
$\varepsilon(2s)=n_{2s-1}$ and $\varepsilon(2s-1)=-n_{2s}$.

\begin{prop}\label{prop:fox}  The Fox derivatives associated to
$\varphi\in {\cal E}(q)$ define a  $\Z [\Z
_2]$--linear map
denoted by
$$d_\varphi\co M_2=\Sigma_{1\leq i\leq 2g}\Z[\Z _2]R_i+\Z[\Z _2]
R_0\to M_1=\Sigma_{1\leq i\leq 2g}\Z[\Z _2]u_i+ \Z[\Z _2]u_0.$$
Its matrix,  with coefficients in $\Z[\Z _2]$, has the
following form:
\begin{equation*}
\left( \begin{array}{cccccll}
 1-t &0&\cdots&0 &0&\varepsilon (1)(1-t) \\
 0 &1-t&\cdots &0 &0&\varepsilon (2)(1-t) \\
 \vdots &\vdots&\vdots &\vdots &\vdots & \vdots\\
0 &0&\cdots&1-t&0&\varepsilon (2g-1)(1-t) \\
0&0& \cdots& 0& 1-t&\varepsilon (2g)(1-t)\\
n_1(1-t)&n_2(1-t)&\cdots &n_{2g-1}(1-t)&n_{2g}(1-t)&c(1+t)
\end{array}
 \right)
\end{equation*}
\end{prop}

\begin{proof}
 The coefficients $m_{ji}$ are:
$$\eqalignbot{
m_{ii}&=\varphi q(1-u_iu_0u_i^{-1})=\varphi (1-u_0)=1-t, i\neq 0;\cr
m_{ji}&=0,i\neq j,i\neq 0,j\neq 0;\cr
m_{0i}&=\varphi q(u_i-[u_i,u_0])=\varphi (u_i)-1,i\neq 0;\cr
m_{(2j -1),0}&=1-\varphi (u_{2j});\cr
 m_{2j,0}&=\varphi (u_{2j -1})-1;\cr
m_{00}&=c(1+t).
}\proved$$
\end{proof}

The relation $\Sigma_{i=1}^{2g} n_i\varepsilon (i)=0$ implies that
the  $\Z [\Z_2]$--module $\mathrm{Im}\ d_\varphi$ is generated by $\{ (1-t)v_i,
1\leq i\leq 2g \}$ and $c(1+t)u_0$ with $v_i=u_i+n_iu_0$.

\begin{nota} ${\bf V}:=\{v_i,1\leq i\leq 2g, v_0=u_0\}$ and ${\bf 
Q}:=\{ R_1,\cdots ,R_{2g},Q\},Q=R_0-\Sigma \varepsilon(i)R_i$.  Also
${\bf{\bar V}}$ is the notation for $\bf{V}$ modulo $\mathrm{Im}\  d_\varphi$.
\end{nota}
 
The structure of  $\Z[\Z _2]$--module of  $H_\varphi$ is:
\begin{align*}
 (H_\varphi,{\bf{\bar{ V}}})= \bigoplus_{1\leq i\leq
2g}&\dfrac{\Z[t]}{((1-t^2),(1-t))}{\bar v_i}\oplus
\dfrac{\Z[t]}{((1-t^2),c(1+t))}{\bar u_0};\\
&\dfrac{\Z[t]}{((1-t^2),(1-t))}\simeq \dfrac{\Z[t]}{(1-t)}\simeq \Z;\ \  \dfrac{\Z[t]}{((1-t^2),c(1-t))}\simeq\dfrac{\Z}{c\Z}\times\Z.\end{align*}

\begin{defi}\label{defi:Dma}
The matrix of $d_\varphi\otimes \mathrm{id}_{\Z_2}\co (M_2\otimes\Z _2,{\bf
R})\to (M_1\otimes\Z _2,{\bf U})$ is called the derived matrix
associated to $\varphi\in {\cal E}(q)$.
\end{defi}

This matrix  is
\begin{equation*}*\label{eqi:Dmatr}
(1+t)\left( \begin{array}{ccccl}
 1 &0&\cdots &\cdots &n_2  \\
 0 &1&\cdots &\cdots &n_1  \\
 \vdots &\vdots &\vdots  &\vdots &\vdots \\
0&0& \cdots & 1 &n_{2g-1}\\
n_1&n_2&\cdots &n_{2g} &c \mod 2
\end{array}
 \right)
\end{equation*}
 with $n_i=\varphi(u_i)\in\{0,1\}$.

\begin{prop}\label{prop:evod}
{\rm (1)}\qua The following sequence is exact:
\begin{equation*}
0\to M_2\otimes\Z _2\stackrel{d_\varphi\otimes \mathrm{id}_{\Z _2}}
\longrightarrow M_1\otimes\Z _2\to H_\varphi\otimes\Z _2\to 0.
\end{equation*}
{\rm (2)}\qua When $c$ is odd, the  matrix of
$d_\varphi\otimes \mathrm{id}_{\Z_2}\co (M_2\otimes\Z _2,{\bf Q})\to
(M_1\otimes\Z _2,{\bf V})$
 is $(1-t) \mathrm{Id}_{2g+1}$,
 and 
$$(H_\varphi\otimes Z_2,{\bf{\bar V}})=
\oplus_{1\leq i\leq 2g}\Z_2{\bar v_i}\oplus\Z_2{\bar u_0}.
$$
{\rm (3)}\qua When $c$ is even, the matrix of 
$d_\varphi\otimes \mathrm{id}_{\Z_2}\co
(M_2\otimes\Z _2,{\bf Q})\to (M_1\otimes\Z _2,{\bf V})$  is
$$(1+t)\left( \begin{array}{cc}
                    I_{2g}&0\\
                    0&0
                    \end{array}\right), $$
$$\displaylines{\hbox{then}\quad\hfill H_\varphi\otimes Z_2\simeq \oplus_{1\leq i\leq 2g}\Z_2{\bar
v_i}\oplus\Z_2[\Z_2]{\bar u_0}.\hfil\qed}$$
\end{prop}

\subsection{Congruence of derived matrices}

Let $\varphi$ and $\varphi'$ be elements of ${\cal E}(q)$, and consider the
following diagram:
\begin{equation}\label{eqi:diagm}
\begin{aligned}
\xymatrix{
(M_2\otimes \Z_2,{\bf R })\ar[d]_{\theta} \ar[rr]^{d_\varphi\otimes
\mathrm{Id}_{\Z_2}}&&(M_1\otimes \Z_2,{\bf U})\ar[d]_{\psi }  \\
(M_2\otimes \Z_2,{\bf R})\ar[rr]^{d_{\varphi'}\otimes
\mathrm{Id}_{\Z_2}}&&(M_1\otimes \Z_2,{\bf U})
}\end{aligned}
\end{equation}
 where the matrix of $\psi$ in the basis ${\bf U}$ is $J(A),
A\in Sp(\Z _2,2g)$ (see \fullref{prop:J} for the definition of $J$).

The $\Z_2[\Z _2]$--map $\theta $   is supposed invertible.
Its matrix is denoted by
$$B=\left(\begin{array}{cc}
 B_1&B_2\\
 B_3&b
\end{array}
 \right),
 \hbox{ with }B_1=(b_{ij};i,j\leq 2g), B_2=\left( \begin{array}{c}
                                                   c_1\\
                                                   \vdots\\
                                                   c_{2g}
                                                   \end{array}
                                                   \right) \hbox{ and }
B_3=(b_1,\cdots,
b_{2g}).$$

We write $n_i:=\varphi(u_i), n'_i:=\varphi'(u_i)$, $\varepsilon(2s):=n_{2s-1}, \varepsilon'(2s):=n'_{2s-1}$ and
$\varepsilon(2s-1):=n_{2s}, \varepsilon'(2s-1):=n'_{2s}$.

\begin{rem}\label{rem:actJ} The condition
``the matrix of $\psi$ in the basis ${\bf U}$ is $J(A), A\in Sp(\Z
_2,2g)$'' implies that  the inverse of $\psi$ in the basis ${\bf U}$ is
also in $J(Sp(\Z _2,2g))$.
\end{rem}

\begin{prop}\label{prop:dn} Let $\psi$ and $\theta$ be as above.
The diagram \eqref{eqi:diagm} is commutative if and only  if the
parameters  verify the following conditions ${\rm mod}\ (1+t)$
\begin{align}
b_{ij}& =a_{ij}+b_j\varepsilon '(i) \tag{\hbox{$\alpha$}}\\
c_{i}&=\Sigma a_{ij}\varepsilon (j)+b\varepsilon'(i)\tag{\hbox{$\beta$}}\\
w_{j}+n_j& =\Sigma n'_ia_{ij}+b_jc \tag{\hbox{$\gamma$}}\\
0&=(1+b+\Sigma b_j\varepsilon (j)),\tag{\hbox{$\delta$}}
\end{align}
with $w_j=\Sigma a_{i,j}r_i+ S_j+r_j$, $S_j=\Sigma a_{2i,j}a_{2i-1,j}$.
\end{prop}

\begin{proof}
 Mod $(1+t)$, the commutativity of the diagram \eqref{eqi:diagm} gives
the following equations:
\begin{align}
a_{ij}& =b_{ij}+b_j\varepsilon '(i)\tag{\hbox{$\alpha$}}\\
\Sigma a_{ij}\varepsilon (j)& =c_{i}+b\varepsilon'(i)\tag{\hbox{$\beta$}}\\
w_{j}+n_j& =\Sigma n'_ib_{ij}+b_jc\tag{\hbox{$\gamma'$}}\\
\Sigma w_{i}\varepsilon (i)+c& =\Sigma n'_jc_{j}+bc.\tag{\hbox{$\delta'$}}
\end{align}
Using the fact that $\Sigma n_i\varepsilon (i)=0,\Sigma n'_i\varepsilon'
(i)=0$, the equation $ (\alpha)$ implies that $\Sigma n'_ib_{ij}=\Sigma
n'_ia_{ij}$. Hence, the equation $(\gamma' )$ is now $(\gamma
):w_{j}+n_j =\Sigma n'_ia_{ij}+b_jc$. The equation $(\gamma )$ implies that
$\Sigma_j w_{j}\varepsilon (j)=\Sigma_{i,j} n'_ia_{ij}\varepsilon
(j)+(\Sigma b_j\varepsilon (j))c$ and $ (\beta )$ implies that
$\Sigma_{i,j} n'_ia_{ij}\varepsilon (j)=\Sigma n'_ic_i$. 
Now $(\delta
')$ becomes $c(1+b+\Sigma b_j\varepsilon (j))=0$.

 If $c$ is odd,  the relation $(\delta)$ is true.

If $c$ is even, we have to add the relation $(\delta):0=(1+b+\Sigma
b_j\varepsilon (j)),\quad  \mod (1+t)$ which is the condition to get the
invertibility of the matrix $B$. This is obtained from the following
computations:

Write $B=B_0+(1+t)K$ with $B_0\in GL(\Z_2,2g)$.  An element $(x_1,\cdots
,x_{2g},x_0)\in \ker B_0$ verifies, $\mod (1+t)$
\begin{equation*}
\begin{array}{rl}
\forall i, \Sigma_j (a_{ij}+b_j\varepsilon '(i))x_j+ (\Sigma_j
a_{ij}\varepsilon (j)+b\varepsilon '(i))x_0&=0\\
\Sigma b_jx_j+bx_0&=0.
 \end{array}
\end{equation*}
The matrix $ (a_{ij})$ is invertible, so for all $j$,
$ x_j=\varepsilon (j)x_0$ and $x_0(\Sigma b_j\varepsilon
(j)+b)=0$. This proves that
$B_0$ is bijective if and only if $b=1+\Sigma b_j\varepsilon (j) \mod
(1+t)$.

Moreover, we have that $B$ bijective if and only if $B_0$ is bijective.
One implication is evident. To prove the converse, let us write
$B=B_0+(1+t)K$ with $B_0\in GL(\Z_2,2g)$, then
$(B_0^{-1}B)^2=(\mathrm{Id}+(1+t)B_0^{-1}K)^2= \mathrm{Id}$, as a matrix with entries in
$\unpfrac{\Z_2[t]}{1-t^2}$. So $B_0^{-1}BB_0^{-1}$ is the inverse of $B$.
\end{proof}

\begin{rem}\label{rem:eo} Once  chosen the basis ${\bf U}$, 
a symplectic matrix
$A=(a_{i,j})$ and any pair $\varphi,\varphi'$, we have:

(1)\qua If $c$ is even, $(\gamma)$ becomes $w_{j}+n_j =\Sigma n'_ia_{ij}$,
which involves a relation between $\varphi$ and $\varphi'$, which is the
necessary and sufficient condition for the existence of $\theta$.

(2)\qua If $c$ is odd, we choose $b_j$ such that $(\gamma)$ is fulfilled, then
$b_{ij}$ and $b$ such that $(\alpha)$ and $(\delta)$ are true and then
$c_i$. This means that it is possible to find an isomorphism $\theta$ such
that the diagram \eqref{eqi:diagm} commutes, hence we have the following
definition:
\end{rem}

\begin{defi}\label{defi:congr}
Let $\varphi$ and $\varphi'$ be elements of ${\cal E}(q)$, the derived matrices
$d_\varphi\otimes \mathrm{Id}_{\Z_2}$ and $d_{\varphi'}\otimes \mathrm{Id}_{\Z_2}$ are
said congruent via $(\psi,\theta)$ if there exist
$\Z_2[\Z_2]$--isomorphisms $\psi$ and $\theta$ such that the following
diagram commutes:
\begin{equation*}
\xymatrix{
(M_2\otimes \Z_2,{\bf R })\ar[d]_{\theta} \ar[rr]^{d_\varphi\otimes
\mathrm{Id}_{\Z_2}}&&(M_1\otimes \Z_2,{\bf U})\ar[d]_{\psi }  \\
(M_2\otimes \Z_2,{\bf R})\ar[rr]^{d_{\varphi'}\otimes
\mathrm{Id}_{\Z_2}}&&(M_1\otimes \Z_2,{\bf U})
}
\end{equation*}
with the constraints that the matrix of $\psi$ in the basis $\bf U$ is
an element of $J(Sp(\Z _2,2g))$ and the matrix of $\theta$
in the basis $\bf R$ is of the following type: $$\left( \begin{array}{cc}
 B_1&B_2\\
0&1
\end{array}
 \right)
 $$
\end{defi}

With this definition, independently of the parity of $c$, the only
condition remaining to get the congruence of the derived matrices is the
condition $(\gamma)$ of  \fullref{prop:dn}.
So we get the main theorem:

\begin{thm}\label{thm:main}
Two special $2$--fold coverings $\varphi$ and $\varphi'$ are $s$--related
(see \fullref{eqn:star})
if and only if the derived matrices associated to $\varphi$ and
$\varphi'$ are congruent.
\end{thm}

\subsubsection[The \ast-product]{The $\ast$--product}\label{Star}

We need to lift the intersection product from $H_1(F_g;\Z_2)$ to
$(H_\varphi\otimes Z_2,{\bf{\bar V}})$.

Replacing $t$ by $1$ gives the description of the projection
$H_1(E_{\varphi},(E_{\varphi})_0;\Z_2)\to H_1(P;\Z_2 )$. Considering a new
basis
$\tau:=\{\tau_i=\nu_i+\varphi(u_i)\nu_0, 1\leq i\leq 2g;\tau_0=\nu_0\}$
 of $H_1(P;\Z_2)$, we define successively
\begin{eqnarray*}
\pi_\varphi\co (H_\varphi\otimes \Z_2,{\bf{\bar V}})&\longrightarrow
& (H_1(P;\Z_2),\tau);\nonumber \\
\Sigma d_i{\bar v_i}+y(t){\bar u_0}&\mapsto &\Sigma d_i\tau_i+y(1)\tau_0,
\end{eqnarray*}
where $y(t)$ is in fact a constant in $\Z_2$ if $c$ is odd and an
element of $\Z_2[\Z _2]$ if $c$ is even,
and $p_{\varphi}=p_\star\circ \pi_\varphi$ the composition of the
projections
$$H_\varphi\otimes\Z_2\stackrel{\pi_\varphi}\longrightarrow
H_1(P;\Z_2)\stackrel{p_\star}\to H_1(F_g;\Z_2).$$

\begin{defi}\label{defi:ast}
A  product, denoted by $\ast$, is defined  in $H_\varphi\otimes\Z_2$ by
lifting the intersection product in $ H_1(F_g;\Z_2)$:
$$x,y\in  H_\varphi\otimes\Z_2\mapsto x\ast y=
p_{\varphi}(x).p_{\varphi}(y)\in \Z_2,$$
where $a.b$ is the intersection product of two elements of $ H_1(F_g;\Z_2)$.
\end{defi}

\subsubsection[Preserving the  ast-product]{Preserving the  $\ast$--product}\label{preStar}
Suppose that $\varphi$ and $\varphi'$ are two special $2$--fold coverings
and $\Psi\co (H_\varphi\otimes\Z_2,{\bf{\bar V}})\to (H_{\varphi
'}\otimes\Z_2,{\bf{\bar V'}})$ is a $\Z_2[\Z_2]$--isomorphism. 
Let us denote by
$
\big(\begin{smallmatrix}
A&B\\
C&D
\end{smallmatrix}\big)
$
the matrix of $\Psi$.
Here
$A$  is a $(2g\times 2g)$--matrix, $B$ is a column with coefficients in
$\unfrac{\Z_2[t]}{(1-t)}\simeq \Z_2$, $C$ is a line and $D$ is an element
in $\unfrac{\Z_2[t]}{((1-t),c(1+t))}$. This ring $\unfrac{\Z_2[t]}{((1-t),c(1+t))}$ is isomorphic to $ \Z_2$ if $c$
is odd, and to $ \unfrac{\Z_2[t]}{(1-t^2)}$ if $c$ is even. 

A generator of $\ker  p_{\varphi}$ is $\tau_0= \nu_0$ and for all $V,
\nu_0\ast V=0$ and ${\bar v_i}\ast {\bar
v_j}=p_\star(\nu_i).p_\star(\nu_j)=\sigma_i.\sigma_j;$ hence we have the
following proposition:

\begin{prop}\label{prop:ast}
A $\Z_2[\Z_2]$--isomorphism $\Psi\co (H_\varphi\otimes\Z_2,{\bf{\bar
V}})\to (H_{\varphi '}\otimes\Z_2,{\bf{\bar V'}})$ respects the
product, ie,
$\Psi (x)\ast \Psi (y)=x\ast y \in \Z_2$,  if and only if  there exists a
symplectic isomorphism $f\co H_1(F_g;\Z_2)\to H_1(F_g;\Z_2)$ such that
$$f\circ  p_{\varphi}= p_{\varphi'}\circ \Psi.$$
\end{prop}

\begin{proof}
$\Psi (x)\ast \Psi (y)=x\ast y$ if and only if $A\in Sp(\Z _2,2g)$ and
$B=0$. These conditions  are equivalent to
$\Psi(\ker p_\varphi)=\ker p_{\varphi'}$ and the existence of such a
symplectic map $f\co H_1(F_g;\Z_2)\to H_1(F_g;\Z_2)$ such that
$$f\circ  p_{\varphi}= p_{\varphi'}\circ \Psi.\proved$$
\end{proof}

Let  $f\co (H_1(F_g;\Z_2),\sigma)\to (H_1(F_g;\Z_2),\sigma)$ be a
symplectic isomorphism with $A=(a_{ij})$ as symplectic
 matrix in the basis $\sigma$.

Let us denote by:
\begin{itemize}
\item $\Psi_f$ the isomorphism from $(H_\varphi\otimes\Z_2,{\bf{\bar V}})$
to $(H_{\varphi '}\otimes\Z_2,{{\bf\bar V'}})$, with matrix
$
\big(\begin{smallmatrix}
A&0\\
0&1
\end{smallmatrix}\big). $ We have $f\circ p_\varphi=p_{\varphi'}\circ \Psi_f;$

\item $\psi_f$ the automorphism of $(M_1\otimes \Z_2, U)$ with matrix
$J(A)$ (see \fullref{defi:tildef} and 
\fullref{prop:J}). Its matrix in the basis $(V,V')$ is $
\big(\begin{smallmatrix}
A&0\\
M&1
\end{smallmatrix}\big), $ with $M=(w_j+n_j+\Sigma a_{i,j}n'_i)$.
\end{itemize}

If $c$ is even, then
\begin{equation}\label{equation**}
\quad \psi_f(\mathrm{Im}\ d_\varphi)=\mathrm{Im}\ d_{\varphi'}\tag{\hbox{$**$}}
\end{equation}

if and only if $M=0$. If so, the quotient isomorphism is equal to
$\Psi_f$ and there exists an isomorphism $\theta$ as in 
\fullref{defi:congr}.

If $c$ is odd, then the relation \eqref{equation**} is always true for any $M$ and
the quotient isomorphism, in the basis $(V,V')$, is
also
$
\big(\begin{smallmatrix}
A&0\\
M&1
\end{smallmatrix}\big). $  Nevertheless $M=0$ is the condition to be added
for getting an isomorphism $\theta$ as in \fullref{defi:congr}.

It is possible to synthesize this study into a definition:

\begin{defi}
Let $\Psi$ be an isomorphism of $H_\varphi\otimes\Z_2$ to $H_{\varphi
'}\otimes\Z_2$ respecting the product,  and $f$ the symplectic
isomorphism of  $H_1(F_g,\Z_2)$ associated by  \mbox{\fullref{prop:ast}}. We will say that $\Psi$ is a quotient if the following
conditions are fulfilled: $\Psi$ is equal to $\Psi_f$ and is a quotient
isomorphism of $\psi_f$. (When $c$ is even, these two conditions are
equivalent).
\end{defi}

\begin{thm}\label{thm:fpsi}
There exists a
 $\Z_2[\Z_2]$--isomorphism
$$\Psi\co
H_1(E_{\varphi},(E_{\varphi})_0;\Z_2)\to
H_1(E_{\varphi'},(E_{\varphi'})_0;\Z_2)$$
 which  is a quotient if and
only if $\varphi$ and $\varphi'$ are $s$--related.
\end{thm}

\subsubsection[Effect of a change of generators of the fundamental group of P on the derived matrices associated to some varphi in E(q)]{Effect of a change of generators of $\pi_1P$ on the
derived matrices associated to some $\varphi\in {\cal E}(q)$}
 The derived matrix associated to $\varphi\in {\cal E}(q)$ is the matrix of the
linear map  $d_\varphi\otimes \mathrm{id}_{\Z_2}\co (M_2\otimes\Z _2,{\bf
R})\to (M_1\otimes\Z _2,{\bf U})$ defined by
$$d(R_j)=\Sigma\varphi\Big(\dfrac{\partial R_j}{\partial u_i}\Big)u_i.$$ Comparing
with \fullref{sec:qo}, we forget the map $q\co \Z_2[F]\to
\Z_2[\pi_1P]$, where $F$ is the free group with $2g+1$ generators.
Suppose that $(u'_i)_{0\leq i\leq n}$ is another choice of generators of
$\pi_1P$ such that for each $i$, $u_i=w_i(u'_0,\cdots ,u'_n)$ is a word.
We are in the situation where, if $R_j=W_j(u_1,\cdots ,u_n,u_0)$, the
new relations are
$R'_j=W_j(w_1,\cdots ,w_n,w_0)$ and $
H_1(E_{\varphi},(E_{\varphi})_0;\Z_2)=\oplus\Z_2[\Z_2]u'_i/ \mathrm{Im}\,
d'_\varphi\otimes \mathrm{id}_{\Z_2}$ with $$d'_{\varphi}(R'_j)=\Sigma\varphi
\Big(\dfrac{\partial R'_j}{\partial u'_i}\Big)u'_i.$$  By induction on the length
of the word $W_j$, it is possible to prove that $$\dfrac{\partial
R'_j}{\partial u'_i}=\Sigma_k\dfrac{\partial R_j}{\partial
u_k}\dfrac{\partial u_k}{\partial u'_i}.$$ Let us denote by $C$ the
matrix with entries $(\unfrac{\partial u_j}{\partial u'_i})$, $M$ and
$M'$ the matrices of
$d_\varphi\otimes \mathrm{id}_{\Z_2}\co (M_2\otimes\Z _2,{\bf R})\to
(M_1\otimes\Z _2,{\bf U})$ and $d'_\varphi\otimes \mathrm{id}_{\Z_2}\co
(M_2\otimes\Z _2,{\bf R'})\to (M_1\otimes\Z _2,{\bf U'})$. We have the
relation:
$$M'=\varphi(C)M.$$
Let us also remark that the matrix $C'$ with entries $\unfrac{\partial
u'_j}{\partial u_i}$ verifies $\varphi(C)\varphi(C')=\mathrm{Id}$ so
$M=\varphi(C')M'$.
Here two systems of generators must be the  lifts of a fixed choice of
generators of $\pi_1F_g$, hence
they are  related by $u'_i=u_0^{-\alpha_i}u_i$ (or equivalently
$u_i=u_0^{\alpha_i}u'_i$), $\alpha_i\in \{0,1\}$.
The matrix $M'$ of $d'_\varphi\otimes \mathrm{id}_{\Z_2}\co (M_2\otimes\Z
_2,{\bf R'})\to (M_1\otimes\Z _2,{\bf U'})$ may be considered as the
matrix of $d_{\varphi'}\otimes \mathrm{id}_{\Z_2}\co (M_2\otimes\Z _2,{\bf
R})\to (M_1\otimes\Z _2,{\bf U})$ with $\varphi'(u_i)=\varphi(u_i)+
\alpha_i, \varphi'(u_0)=\varphi(u_0)$.  The effect is like changing the quadratic
section $s$ (see \fullref{rem:ns}).

The conclusion is that the only invariants,  (independent of the choice
of the generators of $\pi_1P$ , lifting some fixed canonical system of 
generators of $\pi_1F_g$), are the number of classes under the
$s$--relation and the number of elements in each class.

\bibliographystyle{gtart}
\bibliography{link}

\end{document}